\def \version {2014--11--19}

\documentclass[12pt]{article}

\usepackage{amssymb}

\def \DDG {Disjoint Domination Game}
\def \BDG {Bicolored Domination Game}
\def \p {p}
\def \bb {b}
\def \pur {purple}
\def \bbl {blue}
\def \ovc {\overline{c}}
\def \ovb {\overline{\bb}}
\def \ovp {\overline{\p}}
\def \Vp {V_\p}
\def \Vb {V_\bb}
\def \send {$\langle s^* \rangle$}
\def \dend {$\langle d^* \rangle$}
\def \ONS {{\sf ONS}}
\def \ONst {{Opposite Neighbor Strategy}}
\def \ONSP {{\sf OPS}}
\def \ONPst {{Opposite-to-Previous Strategy}}
\def \Do {Dom}
\def \Se {Sepy}

\newtheorem{Theorem}{Theorem}
\newtheorem{lem}[Theorem]{Lemma}
\newtheorem{defi}[Theorem]{Definition}
\newtheorem{crl}[Theorem]{Corollary}
\newtheorem{prp}[Theorem]{Proposition}
\newtheorem{prm}[Theorem]{Problem}
\newtheorem{cnj}[Theorem]{Conjecture}
\newtheorem{rmk}[Theorem]{Remark}
\newtheorem{xmp}[Theorem]{Example}
\newtheorem{obs}[Theorem]{Observation}

\newcommand{\case}[2]{\msk \nin\underline{Case #1:}\quad #2\par\ssk }

\def \bo {\begin{obs} \ }
\def \eo {\end{obs}}

\def \bp {\begin{prp} \ }
\def \ep {\end{prp}}

\def \bpm {\begin{prm} \ }
\def \epm {\end{prm}}

\def \bc {\begin{crl} \ }
\def \ec {\end{crl}}

\def \bcj {\begin{cnj} \ }
\def \ecj {\end{cnj}}

\def \bpm {\begin{prm} \ }
\def \epm {\end{prm}}

\def \thm {\begin{Theorem} \ }
\def \ethm {\end{Theorem}}

\def \bl {\begin{lem} \ }
\def \el {\end{lem}}

\def \bd {\begin{defi} \ \rm }
\def \ed {\end{defi}}

\def \brm {\begin{rmk} \ }
\def \erm {\end{rmk}}

\def \bxm {\begin{xmp} \ \rm }
\def \exm {\end{xmp}}

\def \nmr {\begin{enumerate}}
\def \enmr {\end{enumerate}}

\def \tmz {\begin{itemize}}
\def \etmz {\end{itemize}}

\def \smin {\setminus}

\def \ssq {\subseteq}
\def \nin {\noindent}
\def \bsk {\bigskip}
\def \msk {\medskip}
\def \ssk {\smallskip}
\def \pf {\nin{\bf Proof.} \ }
\def \prf {\nin{\it Proof.} \ }
\def \qed {\hfill $\Box$}

\def \nnn {{\mathbb{N}}}

\def\cC{{\cal C}}

\def \komm {\Huge \sf}

\newtheorem{claim}{Claim}

\def \bcl {\begin{claim} \ }
\def \ecl {\end{claim}}

\def \dia {\hfill $\Diamond$}

\begin{document}

\title{The \DDG ~\thanks{
  ~Research has been supported by the European Union and Hungary
  co-financed
 by the European Social Fund through the project T\'AMOP-4.2.2.C-11/1/KONV-2012-0004 -- National Research Center
 for Development and Market Introduction of Advanced Information and Communication
 Technologies.}}
\author{Csilla Bujt\'as~$^1$\qquad
   \vspace{2ex}
        Zsolt Tuza~$^{1,2}$\\
\normalsize $^1$~Department of Computer Science and Systems
 Technology \\
  \normalsize University of Pannonia \\
\normalsize  Veszpr\'em,
   \vspace{1ex}
     Hungary \\
\normalsize $^2$~Alfr\'ed R\'enyi Institute of Mathematics \\
       \normalsize Hungarian Academy of Sciences \\
\normalsize  Budapest,
   \vspace{1ex}
     Hungary
 }
\date{\footnotesize Latest update  on \version}
\maketitle

\begin{abstract}

 We introduce and study a Maker-Breaker type game in which
 the issue is to create or avoid two disjoint dominating sets
  in graphs without isolated vertices.
 We prove that the maker has a winning strategy on all \emph{connected}
  graphs if the game is started by the breaker.
 This implies the same in the $(2:\!1)$ biased game also in the
  maker-start game.
 It remains open to characterize the maker-win graphs in the
  maker-start non-biased game, and to analyze the $(a:b)$ biased game
  for $(a:b)\neq (2:\!1)$.
  For a more restricted variant of the non-biased game
 we prove that the maker can win on every graph without isolated
 vertices.

 \bsk

\bigskip

\nin \textbf{AMS 2000 Subject Classification:}
 05C57,       
 91A43,       
 91A46,      
 05C69.       

\end{abstract}

\bsk

\section{Introduction}

It is well known that every graph without isolated vertices
 contains two dominating sets which are disjoint.
In this paper we introduce a combinatorial game in which
 one of the two players aims at constructing two disjoint
 dominating sets, while the other player wants to prevent this.
We investigate who has a winning strategy if the graph of the
 game has a certain structural property.

\subsection{Terminology and preliminaries}

 Throughout this paper, we consider finite, simple graphs
  \emph{without isolated vertices}.
For such a graph
 $G=(V,E)$ and for a vertex $v\in V$, the open neighborhood $N(v)$ of $v$ is
 the set of all vertices adjacent to $v$, and its closed
 neighborhood is $N[v]=N(v) \cup \{v\}$. Each vertex dominates
 itself and its neighbors, moreover a set $S \subseteq
 V$ dominates the vertices contained in $N[S]=\bigcup_{v\in S}
 N[v]$. A \emph{dominating set} of $G$ is just a subset $D$ of $V$ which
 dominates all vertices of the graph; that is, $N[D]=V$. The minimum
 cardinality of a dominating set  is termed the domination number of
 $G$.

 Domination is a well-studied subject in graph theory, with many related
 applications. A general overview can be found in \cite{Domination}.
 On the other hand, domination in hypergraphs (set systems) is a
 relatively new area; see \cite{JT09} and \cite{BHT} for results and references.

 \bsk

 Recently, Bre\v{s}ar,  Klav\v{z}ar and   Rall
 \cite{BKR-SIAM} introduced the concept of the \emph{domination game}.
  It is played on a   graph $G$ by  two players, named  Dominator and
 Staller. They  take turns choosing a vertex from $V$ such
 that  at least one new vertex must be dominated in each turn.
 The game ends when no more legal moves can be taken.
 In this game  Dominator's aim is to finish the game with a small dominating set,
  while   Staller aims to
  delay the end of the game.
 The game domination number
  is the number of turns in the game
  when the first turn is Dominator's move  and both  players play
  optimally. For detailed description and  results
  on this subject, see the papers
  \cite{BDKK,BKKR,BKR-SIAM, BKR, JH,CS1,CS2,DKR, KWZ, Kos}.
    Let us mention further that a version of this game for
     the total dominating sets (where $\bigcup_{v\in D} N(v)
 = V$ is required for the \emph{open} neighborhoods)
  was also introduced in \cite{HKR, HKR2}.

  \subsection{\DDG}

 We   define the \DDG\ (DDG, for short) as a two-player game,
  where the players are named Dom and
 Sepy---these are the shortened forms of Dominator and Separator.
 For the game, we have an isolate-free graph $G$ and a color palette
  $\cC= \{ \p, \bb \} = \{ \pur, \ \bbl \}$.
  In the game, Dom and Sepy take turns choosing a vertex and assigning it with a color from
 $\cC$. At any stage of the game,   $\Vp$ and $\Vb$ denote the set of vertices
 colored with $\p$ and
  $\bb$, respectively. Moreover, for $c\in \cC$ we denote by $\ovc$ the complementary
 color for which
   $\{ c, \ovc \} = \cC$ holds. Thus,
 $\ovb = \p$ and $\ovp = \bb$.

 The choice of a vertex $v$ and its coloring with a color $c \in \cC$
  is a legal   (or feasible) move in the game if and only if
 \tmz
 \item[$(i)$] $v$ has not been chosen and assigned with a color up to this
 point, that is $v\notin V_p\cup V_b$; and
 \item[$(ii)$] there exists a vertex $u \in N[v]$ which has not been
 dominated in color $c$, that is    $N[u] \cap V_c=\emptyset$.
 \etmz

Note that each player must select a vertex on his turn whenever a
 legal move is available.
We shall discuss the situation in Section \ref{s:pass} for the game
 where some player may pass.

The game terminates when one of the following two situations
 is reached:
\tmz
\item[\send ]
 some vertex has a monochromatic closed neighborhood
\item[\dend ]
 both $\Vp$ and $\Vb$ are dominating sets
\etmz
 The winner is \Se\ if \send\ is reached, and
 \Do\ wins if \dend\ is reached.

 In other words, the aim of \Do\ is to obtain two disjoint
 dominating sets $V_p$ and $V_b$ at the end of the game, whilst \Se\
 would like to prevent him from reaching this situation.

 It follows from the definition of
  dominating set that both players cannot win
 simultaneously. First of all
  we prove that the game always ends with a win of
 one of them.

\bl   \label{s-or-d}
 As long as neither\/ \send\ nor\/ \dend\ is reached,
  the next player has at least one feasible move.
\el

 \pf
Assume that neither \send\ nor \dend\ has been reached. Then, we
have a vertex $v$ which is
 not dominated by $V_c$ (for some $c\in\cC$) or equivalently,
 $N[v]\cap V_c = \emptyset$. As $N[v]$ is not monochromatic in $\ovc$, there
 exists a vertex $u \in N[v]$ which has not been colored up to this
 moment.
Thus, selecting $u$ and assigning $c$ to it is a feasible move
 because then $v$ becomes dominated in $c$.
 \qed
 \bsk

 In Sections \ref{s:sst} and \ref{s:dst} we study the Sepy-start
  and the Dom-start versions of the
 \DDG. Especially, we prove that for every isolate-free graph $G$, Dom
  has a winning strategy whenever Sepy starts the game and the graph is connected,
 but  in the Dom-start  version it depends on
 the given $G$ which player has a winning strategy.
 We also touch the biased version of the game in the short Section \ref{s:bias}.

 In Section \ref{s:bdg}, we introduce a variant called \BDG\ (or shortly BDG).
 While in the \DDG\ both players are allowed to use
  both colors, in BDG  each player has his private color
  not usable by the other player. For this variant, the definition
  of  \dend\, when \Do\ wins, is slightly modified as follows.
  \tmz
  \item[$\langle d^{**} \rangle$ ]
 For a color $c \in \cC$,  $V_c$ dominates all vertices of $G$, and no vertex $v$ has
 its closed neighborhood entirely  contained in $ V_c$.
 \etmz
 In Section \ref{s:bdg}, we  give an explanation for this change.
   The main result of that section is that \Do\ has a winning strategy
 on all graphs, no matter who starts the \BDG.
 \bsk

We close the Introduction with a lemma which is valid in both the
Disjoint and the Bicolored Domination Game.

\bl \label{no-mono}
 If a vertex\/ $v$ of an isolate-free graph\/ $G$ is chosen  in a turn, then\/
  $N[v]$ cannot be monochromatic after a legal coloring of\/ $v$.
\el
 \pf By the condition $(ii)$, if $v$ is colored with $c$ then at
 least one $u \in N[v]$ has not been dominated by $c$ up to this
 moment. If $u=v$, then no vertex from $N(v)\neq \emptyset$ has color
 $c$;
 and if $u\neq v$, then either $u$ has not been selected yet or
 it has been colored with $\ovc$. In either case, $N[v]$ is not
 monochromatic after the move. \qed
 \bsk

\section{The \Se-start \DDG}   \label{s:sst}

For this game, in this section we prove that
 there exists a winning strategy  for \Do\ which works
on each connected graph.

\thm   \label{S-conn}
 If\/ $G$ is connected, then \Do\ has a winning strategy for the \Se-start \DDG\ on\/ $G$.
\ethm

 \pf
After each move of \Se, Dom  applies the following strategy.

\bsk

 \nin {\bf \ONst\ (\ONS )}
  \tmz
  \item
 Suppose that \Se\ selected a vertex $v$ and colored it with $c$.
 Then Dom chooses a vertex according to the following rules.
 \item[]\hspace{-1em}$(\ONS1)$\hspace{1em}\Do\ selects a neighbor
  $u$ of $v$ which can be colored
 with $\ovc$.
 \item[]\hspace{-1em}$(\ONS2)$\hspace{1em}If $(\ONS1)$ cannot be applied, \Do\ selects a
 vertex $u$ which can be colored with $c'\in \cC$, moreover $u$ has
 a neighbor of color $ \overline{c'}$.
\etmz
  First, we observe the following consequences of the definition
  above.

  \bl \label{ONS-opposite}
  As  long as \Do\  does not violate\/ \ONS,
  \tmz
  \item[$(i)$]
  after each move of \Do, every colored vertex has a neighbor assigned
with the complementary color;
\item[$(ii)$] the game cannot terminate with\/ \send.
\etmz
  \el
  \prf
  Remark that if rule   \ONS1 cannot be applied,
      then the vertex
  $v$, colored with $c$ by \Se\ in the last step, already has a
  neighbor of color $\ovc$.
  Then, part $(i)$ of the lemma  immediately follows from the definition  of \ONS.
 Particularly,  no choice of \Do\ can make  any
closed neighborhood $N[x]$  monochromatic. On the other hand, by
Lemma \ref{no-mono} and part $(i)$, when \Se\ selects a vertex $v$,
this cannot result in a monochromatic closed neighborhood either.
Thus, \send\ cannot be reached unless some move of Dom violates
\ONS.
 \dia

\bsk

 Therefore, it suffices to prove that \Do\ can apply \ONS\ in each
 turn.

\bl   \label{dom-go}
 After each step of \Se, \Do\ can apply\/ \ONS\
  as long as neither\/ \send\ nor\/ \dend\ is reached.
\el

 \prf
 We will prove that   there exists a move complying with   \ONS2.
 Let $D_0$, $D_1$, and $D_{\ge 1}$ denote the set of vertices which
 are dominated by 0, precisely 1, and at least 1 color from $\cC$,
 respectively.
 Since \dend\ has not been reached, $D_0 \cup D_1 \neq \emptyset$;
  and $D_{\ge 1} \neq \emptyset$ holds already after the very first move.

 First consider the case $D_0 \neq \emptyset$. As $G$ is connected,
 there exists a vertex $u \in D_{\ge 1}$ with a neighbor $u' \in
 D_0$. Thus, $u$ must be an uncolored vertex which has a neighbor
 $u''$ colored with a $c \in \cC$. Then, choosing $u$ and coloring
 it with $\ovc$ is a legal move for Dom (this dominates   $u'$  with color $\ovc$) and
 also corresponds to   \ONS2.

 In the other case we have $D_0 = \emptyset$, which implies $D_1
 \neq \emptyset$. Hence, there is a vertex $u$ dominated with a
 color $c$, but not dominated with $\ovc$. If $u \notin
 V_c$, it is uncolored and Dom may choose $u$ and assign it with color $\ovc$.
  This satisfies the requirements.
 Now, suppose that $u \in V_c$. As   $u$ is not dominated by $\ovc$ and
 \send\ is not reached, $u$ has an uncolored neighbor $u'$.
 Selecting $u'$ and coloring  it with $\ovc$ is a legal move
 complying with   \ONS2.
 \dia

\bsk

 Lemma \ref{ONS-opposite} and Lemma \ref{dom-go} together mean that
 the game terminates with \dend\  if \Do\ applies the \ONst.
 \qed

\section{The \Do-start game}   \label{s:dst}

Contrary to the \Se-start game, some graphs admit a winning
 strategy for \Se\ in the \Do-start game.

 \bp \label{cycle}
 For every\/ $n \ge 8$, Sepy has a winning strategy for the
 \Do-start \DDG\ on\/ $C_n$.
 \ep

 \pf Let $v_1, v_2,\dots, v_n$ be  the vertices in cyclic order.
 Without loss of generality, we assume that Dom selects $v_1$ and
 assigns it with purple. Then, Sepy can color $v_2$ with purple.
 If Dom's next choice is $v_3$, $v_4$ or $v_5$ (assigning it with either
 color), then Sepy   colors $v_n$ with purple. This is a legal
 choice as $n\ge 8$ ensures that $v_{n-1}$ was not dominated  (with
 any color) before this turn. Thus, $N[v_1]$ becomes monochromatic and
 the game terminates with
 \send.
 Similarly, \send\ is reached if Dom selects a vertex different
  from  $v_3$, $v_4$ and
 $v_5$. In this case, Sepy's next move is coloring $v_3$ with
 purple.
 \qed

 \bp
Let\/ $G$ be a graph of minimum degree at least 2, and let\/ $G^{+2}$
 be the graph obtained by subdividing each edge of\/ $G$ into a
 path of length~3.
Then \Se\ can win the \Do-start \DDG\ on\/ $G^{+2}$.
 \ep

 \pf
We show that these graphs offer \Se\ a local strategy
 to reach \send\ after a small number of steps.
There are two kinds of first moves for \Do :
 to select an original vertex of $G$ or a subdivision vertex.

\case{1}{The first move is a subdivision vertex.}

 \nin
Say, an edge $wz$ of $G$ has been subdivided to a path $wxyz$,
 and \Do\ assigned color $c$ to vertex $x$ in the first move.
Then assigning $c$ to $y$, \Se\ creates a double threat:
 putting color $c$ further on any of $w$ and $z$ terminates
 the game with \send.
\Do\ may color only one of $w$ and $z$ in one step, therefore
 the only way to delay the winning of \Se\ would be to make
 $c$ infeasible on both $w$ and $z$.
This situation would occur precisely if the entire
 $N(w)\cup N(z)$ became dominated by $V_c$ after the second
 move of \Do.
But the neighbors of $w$ and $z$ on the subdivided edges
 do not have any common neighbors, therefore $c$ remains
 feasible for $w$ or $z$ after the move of \Do.

\case{2}{The first move is an original vertex.}

 \nin
Suppose that \Do\ first selects a vertex $w$,
 whose neighbors in $G$ are $z_1,\dots,z_d$, and
 assigned color $c$ to it.
Denote by $wx_iy_iz_i$ the subdivision path
 of the edge $wz_i$ ($i=1,\dots,d$).
The strategy of \Se\ is to color all the vertices
 $x_i$ with $c$ one by one.
Each such move creates a threat on $x_i$, because
 \send\ will occur once \Se\ assigns $c$ to $y_i$.
The only way for \Do\ to prevent this is to color a
 vertex in $N[z_i]$.
 As these sets $N[z_i]$ are pairwise disjoint,
 \Do\ needs a distinct move for each, and hence
 \Se\ can make the entire $N[w]$ monochromatic
 and achieve \send.
 \qed

\bsk

On the other hand, in every DDG played on a complete graph $K_n$ or
on a path $P_n$ ($n \ge 2$)  Dom can win, even if he begins the
game. These examples are special cases of the following theorem.

\thm
 Let\/ $G=(V,E)$ be a connected graph,  which has two different
 vertices\/ $u$ and\/ $v$ satisfying\/ $N[u] \subseteq N[v]$.
 Then, Dom has a winning strategy for the
 \Do-start \DDG\ played on\/ $G$.
 \ethm

\pf (sketch) We say that a vertex $v$ is safe (concerning
DDG) if no matter how the game is continued, $N[v]$ cannot be
entirely monochromatic.

The winning strategy of Dom is as follows.
   \tmz
   \item In the first turn, Dom selects a vertex $v$ which has a neighbor
   $u$ with $N[u] \subseteq N[v]$. If $v$ gets color $c$, then $u$
   cannot be colored with $c$ in any later turns. Hence, $N[v]$
   cannot become monochromatic, this is a safe vertex.
   \item In the later turns, Dom applies \ONS.
  \etmz

  The details of the proof are similar to those of Lemma
   \ref{ONS-opposite} and Lemma \ref{dom-go}.
     After each turn of Dom,  every
   vertex in $V_b \cup V_p$ is safe. Moreover, no move of Sepy
   creates a monochromatic closed neighborhood.
   If \ONS\ has not been violated and the game has not ended,
   Dom can apply \ONS\ in his next turn.
   \qed

\section{Passing allowed}   \label{s:pass}

In the standard versions of domination games both players have to
 move in each turn until no legal move is possible.
It is well known, however, that in some games there exist situations
 where it really does matter if the next player is not allowed to
 skip the move.

By definition we exclude the possibility of passing in the very first
 move, because it would immediately change the character of the game
 by switching between \Do-start and \Se-start.
Another reason for this restriction is that if passing was
 allowed for both players at any time, then in a \Se-win graph in the
 \Do-start game first \Do\ should pass to avoid losing, but then
 also \Se\ should pass because otherwise \Do\ can surely win, as we shall
 see below; hence the game would end up with a trivial draw.
 There exist further situations, too, where the possibility of
 double passing would cause unwanted anomalies.
For this reason we restrict our attention
 to games in which just one specified player --- or none of the
 players --- is allowed to pass.

 Passing may or may not help a player.
Namely, we shall see in Section \ref{s:bdg} that the possibility of passing
 has no effect on the outcome of the \BDG.
  On the other hand, `virtual passing'  may be a useful
 concept for designing strategies in `biased games'
 discussed in the next Section \ref{s:bias}.
Moreover, by comparing the results of the previous two sections
 we can see that passing may help a lot for \Do, as expressed in the
 following variant of Theorem \ref{S-conn}.

\bsk

 \nin
{\bf Theorem \ref{S-conn}$\boldmath '$ }
 \emph{If \Do\ is allowed to pass but \Se\ isn't,
   then \Do\ has a winning strategy in the \Se-start
  \DDG\ on every graph.
  Also, \Do\ has a winning strategy in the \Do-start
  \DDG\ on every graph containing at least one
  \Do-win component}

\bsk

 \pf
For connected graphs the assertion clearly is valid
 by Theorem \ref{S-conn}.
If $G$ is disconnected, then \Do\ can win with the following strategy.
 \tmz
  \item If \Do\ should start the game, then start with the first move of a
   winning strategy in a \Do-win component.
  \item Afterwards, or if \Se\ starts the game, always play in the same
   connected component where \Se\ made his latest move.
  \item In that component apply the \ONst.
  \item If \Se\ moved and there is no legal move in that component
   anymore, then pass.
 \etmz
With the first and the last rules of this strategy \Do\ can force \Se\
 to open each connected component which would be \Se-win in the \Do-start game.
 Consequently, Dom can win  by applying \ONS\
 according to Theorem \ref{S-conn}.
  \qed

\bsk

On the other hand, it turns out that passing has no
 benefit for \Se\   if he begins the game.

\bsk

 \nin
{\bf Theorem \ref{S-conn}$\boldmath ''$ }
 \emph{If\/ $G$ is connected, then \Do\ has a winning strategy for the \Se-start
  \DDG\ on\/ $G$, even when \Se\ is allowed to pass at any time except in the
  first move.}

\bsk

 \pf
\Do\ can win in a similar way as in the original game when passing was
 not allowed.
In fact he never needs to pass, as shown by the
 following scheme.

\bsk

 \nin {\bf \ONPst\ (\ONSP )}
  \tmz
  \item
 Suppose that in the   latest turn a vertex $v$ was selected and colored
   with $c$. Then, Dom chooses a vertex
 according to the following rules.
 \item[]\hspace{-1em}$(\ONSP1)$\hspace{1em}\Do\ selects a neighbor
  $u$ of $v$ which can be colored
 with $\ovc$.
 \item[]\hspace{-1em}$(\ONSP2)$\hspace{1em}If $(\ONSP1)$ cannot be applied, \Do\ selects a
 vertex $u$ which can be colored with $c'\in \cC$, moreover $u$ has
 a neighbor of color $ \overline{c'}$.
 \etmz
It can be verified along the lines of the proof of Theorem \ref{S-conn}
 that this strategy is feasible and \Do\ wins if he applies it.
\qed

\bsk

\section{Biased games}   \label{s:bias}

In a \emph{biased} or \emph{asymmetric} game the players may make
 more than one move at a time.
Such games are parameterized with two positive integers
 indicating the numbers of moves of the players per turn.

Adopting this notion to the \DDG, let $d,s$ be two fixed
 positive integers.
We use the shorthand \emph{$(d:s)$-game} for the game where \Do\
 sequentially selects and colors exactly $d$ vertices in each turn
 (except near the end of the game when only fewer than $d$ possibilities
 remain), and \Se\ colors at most $s$ vertices per turn
  (and may pass if he wishes to do so).
Note that $d$ always refers to \Do\ and $s$ always refers to \Se,
 no matter who starts the game.
The general requirement to dominate new vertices in the colors of the
 successively selected vertices is kept also in the $(d:s)$-game.
Hence the case $d=s=1$ precisely means the \DDG\  where Sepy is
allowed to pass.

Here we only consider the case $s=1$;
 some short comments on larger values will be given
  in the concluding section.
 First, we prove the following lemma:

 \bl \label{safe-comp}
  If Dom has the possibility  to achieve \dend\ by playing just one next vertex
  inside a component
 $ C_i$ then no sequence of legal moves  results in a
 monochromatic $N[v]$ inside $C_i$.
 \el
 \pf
  Assume that the assignment of  color $c$ to $u \in V(C_i)$ would make
  the entire $V(C_i)$   dominated by both colors. Then, already
  without this action, each $w \in V(C_i)$    is dominated by
  $\ovc $. That is, $N[w] \cap V_{\ovc} \neq \emptyset$ and no matter
  which vertex or vertices are chosen (and assigned to c) later, no $N[w]$ will be
  monochromatic in color $c$ in this component.
  On the other hand, in the continuation of the game, no $v \in
  V(C_i)$ can be colored with $\ovc$, which implies that no closed
  neighborhood can become monochromatic in $\ovc$ inside $C_i$.
  \qed
  \bsk

 The components in which no closed neighborhoods can become
 monochromatic in any continuation of the game will be called {\it
 safe components}.

Our main result in this section states that any $d>1$ yields
substantial advantage for \Do.

\thm
  \Do\ has a winning strategy in the\/ $(d:1)$-game
  on every graph, for every\/ $d\ge 2$.
\ethm

 \pf
We advise \Do\ to play a variant of the \ONPst\, with a
 simple modification.
 \nmr
 \item If Dom starts or continues the game on a new component with at least two consecutive choices,
  color  the
 first (any) vertex arbitrarily, and
 then    color the next vertex  according to \ONSP.
 \item
 In the other cases, apply  \ONSP\ itself, apart from one exceptional
  situation described below in the third rule of the strategy.
  Note that when \Do\ applies $\ONSP2$, this current second rule might mean a choice
  of a vertex from another component, from which some vertices were played
  earlier.
  \item If these rules  yield a situation where some components
  $C_1, \dots ,C_i$ are completely dominated in both colors and the
  remaining ones $C_{i+1}, \dots ,$ $C_k$ are completely
  undominated,
  moreover  Dom can take only one   further choice before the turn of Sepy,
  then instead of taking the preceding move in   $C_1\cup\cdots\cup C_i$ he  chooses two
  vertices from the new component $C_{i+1}$ in the way described  by the first rule above.
    By Lemma~\ref{safe-comp},
   the non-completed component inside $C_1\cup\cdots\cup C_i$ is safe.
  \enmr
By the proof of Lemma \ref{dom-go}, this strategy is
 feasible and leads to the winning of \Do.
 \qed

\bsk

\section{The \BDG }   \label{s:bdg}

A natural variant of the \DDG\ is when \Do\ and \Se\ have their
 private colors; that is, \Do\ may only use $\p$ and \Se\ may
 only use $\bb$.
Also in this case the game may terminate with \send\ as above,
 in which case \Se\ wins.
However, recall that we had to modify the meaning of \dend\
slightly:
  \tmz
  \item[$\langle d^{**} \rangle$ ]
 For a color $c \in \cC$,  $V_c$ dominates all vertices of $G$, and no vertex $v$ has
 its closed neighborhood entirely  contained in $ V_c$.
 \etmz
 The point is that if $V_c$ dominates $G$, then the player of color $c$
 does not have any further feasible moves, while the other player may
 still have some; but on the other hand, forbidding $N[v]\ssq V_c$
 for all $v\in V$, the set $V\smin V_c$ dominates $G$ and therefore
 letting the player of color $\ovc$ play as long as a feasible move
 is available, eventually two disjoint dominating sets are reached.
 \def \Mst {Matching Strategy}
\def \MS {{\sf MS}}
\bsk

  In this section we will show that Dom can win the \BDG\ on any isolate-free graph $G$.
 First, we prove this statement for the special case where $G$ contains a perfect
 matching. Then, the \Mst\ introduced here will be extended to
 obtain a winning strategy for Dom which works on   all
 isolate-free graphs.

 \bp   \label{perfmatch}
 If\/ $G$ has a perfect matching, then \Do\ can win
 the Sepy-start and also the \Do-start \BDG\ on\/ $G$.
 \ep

\pf Suppose that $G$ has $2n$ vertices, and let
 $\{v_1v_1',v_2v_2', \dots, v_nv_n'\}$ be a perfect matching in $G$.
The winning strategy of Dom is as follows.

\msk

{\bf \Mst } (\MS )
   \tmz
   \item If Sepy selects a vertex from the matching edge $\{v_i,v_i'\}$,
    say he plays $v_i$, and playing the other vertex $v_i'$ of
    this edge is a legal move for \Do, then \Do\ plays $v_i'$.
   \item Otherwise, Dom is free to make any legal choice
    from any pair $\{v_j,v_j'\}$ in which both vertices
    remained uncolored until that move.
   \etmz
   It is clear that Dom can apply this strategy throughout the game,
   without making any $\{v_i,v_i'\}$ monochromatic at any time.
To show that it is a winning strategy indeed, observe that the only
reason why
 it is not legal for \Do\ to play $v_i'$ as a response to \Se's $v_i$
  is that both $v_i$ and $v_i'$ are already dominated in \Do's color.
At the end of the game this property holds for all edges of the
 perfect matching.
Since \Do\ selects at most one vertex from each matching edge, the
game
 can never end with $\langle s^{*} \rangle$
 but only with $\langle d^{**} \rangle$.
  \qed

\thm
 \Do\ can win the \BDG\ on every graph without
  isolated vertices, both in the \Do-start and \Se-start cases.
\ethm

 \pf  As a preprocessing for his strategy, Dom determines a maximum
 matching in $G=(V,E)$, denote it by $M=\{u_1v_1,u_2v_2,\dots,u_mv_m\}$.
We refer to the edges $e_i=u_iv_i$ as \emph{matching-edges}, and
 their vertices as \emph{matching-vertices}.
In the complementary part $V\setminus M$ we call the vertices
\emph{external}.
 Note that all neighbors of an external vertex are matching-vertices.

Each external vertex is adjacent to some matching-vertex,
 because isolates are excluded.
If some $x\in V\smin M$ is adjacent to both $u_i$ and $v_i$,
 then no other $x'\in V\smin M$ is adjacent to any of $u_i$ and $v_i$,
 by the maximality of $M$.
So, such an $e_i$ is a \emph{triangle-edge}.

On the other hand, if $x\in V\smin M$ is adjacent to $u_i$ but
 not to $v_i$, then \emph{all} $x'\in V\smin M$ having a neighbor
 in $e_i$ have precisely the same neighbor, $u_i$.
  We shall then call $u_i$ a \emph{center}, and $e_i$ a \emph{star-edge}.
Note that if a matching-edge has at least one external neighbor
 then either it is a triangle-edge or its center is uniquely
 determined.

To win on $G$, in both the Dom-start and Sepy-start version,
 Dom applies the following rules.
 \begin{enumerate}
   \item \textbf{Center Rule.}

    If Dom selects a vertex from a star-edge, from which no vertex
     has previously been selected, then he selects the center.

   \item \textbf{Incompleteness Rule.}

    Dom does not select more than one vertex from a matching-edge.

   \item \textbf{Neighbor Rule.}

    After each step of Dom, $V_p$ dominates $V_b$.
        Moreover, if Sepy selected a vertex of some matching edge $e_i$ and
     the other vertex of  $e_i$ is feasible for Dom, then
     Dom plays that vertex. As a consequence, each blue vertex
     selected by Sepy has at least one purple neighbor selected by
     Dom.

   \item \textbf{Matching Rule.}

    Dom selects a vertex outside $M$ only if every feasible
     matching-vertex violates some of the rules above.

 \end{enumerate}

\bl
 If Dom can keep all these rules during the whole game,
  then he wins.
\el

 \pf
We have to prove that the game cannot terminate with \send,
 i.e.\ no closed neighborhood can become monochromatic.

It is clear by the Neighbor  Rule and Lemma~\ref{no-mono}
 that Sepy cannot create
 any $N[v]$ in blue.
Suppose that Dom is forced to do so, say an entire $N[v]$
 becomes purple when Dom selects a vertex $x$.
Again
 by Lemma~\ref{no-mono}, $x\ne v$ certainly holds.
Now, the vertex $v$ with totally purple neighborhood cannot be
 a matching-vertex because its pair in $M$
 cannot be purple together with $v$, by the Incompleteness Rule.
And it cannot be an external vertex either, because it would
 assume that Dom selected the external vertex $v$ earlier
 than the feasible matching-vertex $x$, which would violate
 the Matching Rule.
 \dia

\bsk

As the reader can observe, the Center Rule has not been used
 in the argument of the previous proof.
In fact it is needed to ensure that Dom can keep all the
 other rules during the whole game.

\bl
 As long as\/ $\langle d^{**} \rangle$ has not been reached,
  after any move (or pass) of Sepy, Dom has a legal move
  respecting the four rules above.
\el

 \pf
 Consider a turn of Dom and assume that he obeyed all the four
rules in all of his previous turns. We check the four requirements
  for Dom's current move one by one.

\msk

 \nin \underline{Center Rule.}\quad
Let $e_i=u_iv_i$ be a star-edge with center $u_i$, and suppose that
 none of its ends has been selected yet.
If an external neighbor of $e_i$ is not dominated in purple, then
 Dom can play $u_i$.
  The same move is legal if $u_i$ or $v_i$ is not dominated in purple.
The only bad case is if $V_p$ already dominates $N[u_i]$, but $v_i$
 still has a neighbor, say $x$, which is not dominated in purple.
Since $v_i$ is not the center of $e_i$, this $x$ must be a
 matching-vertex.
It also means that no purple vertex has been selected from
 the matching-edge of $x$.
But then Dom can select the center of that edge
 as a legal move  if it is a star-edge, or any of $x$ and its neighbor
 if it is a triangle-edge.
Note that both vertices of the edge in question were previously
 non-selected, because if Sepy had selected the pair of $x$ in an
 earlier move then Dom would have immediately responded with
 selecting $x$ as a vertex non-dominated in purple, by the
 Neighbor Rule.

\msk

 \nin \underline{Incompleteness Rule.}\quad
This rule would be violated only if Dom selects both
 $u_i$ and $v_i$.
Due to the Center Rule, this can happen only if $v_i$ is adjacent
 to a matching-vertex $y$ non-dominated in purple.
But then the matching-edge containing $y$ either is still completely
 non-selected and Dom can play its center, or the pair of $y$ in
 that edge was selected by Sepy, which contradicts the Neighbor Rule.

\msk

 \nin \underline{Neighbor Rule.}\quad
Suppose first that Sepy selected a vertex of $e_i$ in his
 latest move.
If the entire $e_i$ is not yet dominated in purple, then
 the other vertex of $e_i$ is feasible for Dom, and the
 rule is kept by choosing that vertex.
In the other case the selection of Sepy is already dominated
 in purple, and the rule does not put any restriction on Dom.

 Suppose next that Sepy selected an external vertex $x$.
If $x$ is not dominated in purple yet, then it must have an
 uncolored neighbor, which necessarily is either the center
 of a star-edge or belongs to a triangle-edge.
Thus, Dom can keep the rules.

\msk

 \nin \underline{Matching Rule.}\quad
  This rule is easy to keep.
Dom is forced to select an external vertex only if all neighbors
 of this vertex were already selected by Sepy.
But then Dom can postpone the selection of all such external
 vertices to the last part of the game.
 \dia

\bsk

The two lemmas above together imply that Dom can win the \BDG\
 on every isolate-free graph. \qed

\section{Conclusion}   \label{s:concl}

In this paper we introduced two games on graphs,
 the \DDG\ and the \BDG.
This new area offers many challenging open questions;
 below we collect some of them.

\subsection{\Do-win and \Se-win graphs}

We have proved that \Do\ has a winning strategy in the
 Sepy-start \DDG\ on connected graphs.
It is not very well understood, however, which properties of $G$
 ensure a winning strategy for \Do\ or \Se\ if \Do\ starts
 the game, or if Sepy starts but the graph is disconnected.

\bpm
 Characterize the disconnected graphs on which \Do\ can win the
  \Se-start \DDG.
\epm

\bpm
 Characterize the graphs (connected or otherwise) on which \Do\ can win the
  \Do-start \DDG.
\epm

A false intuition says that,
 due to the constructive goal of \Do,
 the possibility of passing does not increase \Se's chances
 to win.
But in fact this is not true, as shown by the following observation.

\bp
 There exist graphs such that the possibility of passing does
    increase \Se's chances
 to win.
\ep

\pf
  As an example, consider the graph $G$ which is the disjoint
 union of a $C_4$ and a $C_8$. Note that the \DDG\ played on $C_4$
 surely ends with two purple and two blue vertices; no monochromatic
 $N[v]$ can arise, independently of the strategies of the players.
 In addition,  by Proposition~\ref{cycle} and Theorem~\ref{dom-go},
  Sepy can win the Dom-start
 game and Dom can win the Sepy-start game on
 $C_8$ when passing is not allowed.

 Thus, if Sepy starts the game on $G$ and passing is not allowed,
 Dom has a winning strategy as
 he can ensure that Sepy selects the first vertex from $C_8$.
 In contrary, if Sepy may pass, first he can select a vertex from $C_4$
 and then he passes in every turn   until Dom colors a vertex from
 $C_8$. Then, Sepy can win the game.
\qed

\subsection{Biased games}

If $d>1$ or $s>1$, weaker or stronger conditions may be
 imposed than those in the $(d:s)$-game.
Namely, we may allow \Do\ to select fewer than $d$ vertices per turn,
 and/or force \Se\ to select exactly $s$ vertices per turn.
More generally, instead of $d$ and $s$ one may specify
 parameters $d''\ge d'\ge 1$ and $s''\ge s'\ge 0$,
 requiring that, in each turn,
  the number of vertices selected by \Do\ has to be
 between $d'$ and $d''$ while
  the number of vertices selected by \Se\ has to be
 between $s'$ and $s''$.
Even more generally one may specify the sets $D^*$ and $S^*$
 of allowed numbers of selections per turn (possibly varying
 turn by turn), etc.

It is not our goal to analyze the similarities and differences
 between these variants; we leave this direction open for
 future research of other authors.
 It also remains unexplored, which kinds of substructures and
 legal moves should be excluded in order to make the
 $(d:s)$-game non-trivial on some classes of graphs.
  In this direction the following related question seems to be important.

\bpm
 Determine the (sets of) restrictions that ensure the following:
 For every\/ $s\ge 1$ there exists a threshold value\/ $d_s$
  such that \Do\ wins both the \Do-start and \Se-start\/
   $(d:s)$-game on every graph of minimum degree at least\/ $s$.
\epm

Minimum degrees slightly larger than $s$ may also be
 of interest.
Our impression is that one
 natural kind of conditions may be something like this:
 If \Se\ selects a set $S=\{v_1,\dots,v_s\}$
   in a move, then there must exist other $s$ vertices
   $v'_1,\dots,v'_s\notin S$ which have not been
   dominated previously in \Se's color and all of
   $v_1v'_1,\dots,v_sv'_s$ are edges in $G$.
The strategy \ONS\ or some variants of it may turn out to be
 powerful also in this context.


 The following question seems to be interesting, too.

\bpm
 Let\/ $G$ be a graph, and\/ $d,s\in\nnn$.
   Investigate the relation of the\/ $(d:s)$-game
   to the\/ $(d+1:s)$-game and to the\/ $(d:s+1)$-game
 in both the \Do-start and \Se-start versions.
\epm


\subsection{More than two colors}

Some graphs contain more than two mutually disjoint
 dominating sets.
The \emph{domatic number} of $G=(V,E)$ is the largest integer $k$
 for which there exists a vertex partition $D_1\cup\cdots\cup D_k=V$
 into $k$ dominating sets $D_i$ of $G$.

\bpm
 Let\/ $G$ be a graph with domatic number\/ $k$, and let\/ $\cC_\ell$
  be a palette of\/ $\ell\le k$ colors.
 What kind of structural properties of\/ $G$ imply that \Do\ has a
  winning strategy in the game where he and \Se\ alternately assign
  colors from\/ $\cC_\ell$ to the vertices of\/ $G$ and \Do's goal is
  to create a domatic partition with\/ $\ell$ vertex classes?
 In particular, characterize the graphs which admit a winning strategy
  for \Do\ in the case\/ $k=\ell$.
\epm

\end{document}